\documentclass[12pt]{article}

\usepackage{amsmath}
\usepackage{amssymb}
\usepackage{eucal}

\usepackage[english]{babel}

\begin{document}
 
\baselineskip 16pt
 
\title{FINITE GROUPS WITH HALL SCHMIDT SUBGROUPS}


\author{V.\,N. Kniahina and V.\,S. Monakhov}


\maketitle

\begin{abstract}

A Schmidt group is a non-nilpotent group whose every proper
subgroup is nilpotent.
We study the properties of a non-nilpotent 
group $G$ in which every Schmidt subgroup is a Hall subgroup of $G$.

\end{abstract}

{\small {\bf Keywords}: Hall subgroup, Schmidt subgroup, 
non-nilpotent group.}

MSC2010 20D20, 20E34


{\small }

\section{Introduction}

A non-nilpotent finite group whose proper subgroups are 
all nilpotent is called a Schmidt group.
O.\,Yu. Schmidt pioneered the study of such groups \cite{Sm}. 
In a series of Chunihin's papers, Schmidt groups were 
applied in order to find criterions of nilpotency and 
generalized nilpotency, and also to find non-nilpotent 
subgroups, (see \cite{Ch}). A whole paragraph from Huppert's monography 
is dedicated to Schmidt groups, (see \cite{Hup}, III.5). 
Review of the results 
on Schmidt groups and perspectives of its application 
in a group theory as of 2001 are provided in paper \cite{MonUMK}.

Let $S$ be a Schmidt group. Then the following properties hold: 
$S$ contains a normal 
Sylow subgroup $N$ such that $S/N$ is a primary cyclic subgroup; 
the derived subgroup of $S$ is nilpotent; the derived length of  $S$ 
does not exceed $3$; any non-normal Sylow subgroup $Q$ of $S$ is cyclic 
and every maximal subgroup of $Q$ is contained in $Z(S)$;
every normal primary subgroup of $S$ other than a Sylow subgroup of $S$  
is contained in $Z(S)$. 

In this paper the properties of a non-nilpotent 
group $G$ in which every Schmidt subgroup is Hall in $G$
are studied. In particular, for such groups  a number of properties 
of Schmidt groups are applicable. We prove the following theorem.

\medskip

{\bf Theorem.} {\sl Let $G$ be a finite non-nilpotent group
in which every Schmidt subgroup is Hall in $G$.
Then the following statements hold:

1) if $P$ is a non-normal Sylow $p$-subgroup of $G$, then
$P$ is cyclic and every maximal subgroup of $P$ is contained in $Z(G)$;

2) if $P$ is a normal Sylow $p$-subgroup of $G$ and $G$ is 
not $p$-decomposable, then $P$ is either minimal normal in $G$  
or nonabelian, $Z(P)=P^{\prime}=\Phi (P)$,
and $P/\Phi (P)$ is minimal normal in $G/\Phi (P)$;

3) if $P_1$ is a normal $p$-subgroup of $G$, $P_1$ other than a 
Sylow one and $G$ is not $p$-decomposable, then $P_1$ is
contained in $Z(G)$;

4) if $Z(G)=1$, then $G$ has a normal abelian Hall subgroup 
$A$ in which every Sylow subgroup is minimal normal 
in $G$, $G/A$ is cyclic and $|G/A|$ is squarefree.}

\medskip

{\bf Corollary. } {\sl Let $G$ be a finite non-nilpotent group
in which every Schmidt subgroup is a Hall subgroup of $G$.
Then $G$ contains a nilpotent Hall subgroup $H$ such that 
$G/H$ is cyclic. In particular, $G/\Phi (G)$ is metabelian.}


\section{Notations and preliminary results}


Throughout this article, all groups are finite. 
We use the standart terminology
and notations of \cite{Hup}, \cite{Mon}. Recall that a $p$-closed
group is a group with a normal
Sylow $p$-subgroup and a $p$-nilpotent group is a group of order
$p^am$, where $p$ does not divide $m$, with a normal subgroup of order
$m$. A group is called
$p$-decomposable if it is $p$-closed and $p$-nilpotent simultaneously.
A group whose order is divisible by a prime $p$ is a $pd$-group.
We denote by $Z(G)$, $G^{\prime}$, $\Phi(G)$, $F(G)$, $G_p$ the center, 
the derived subgroup, the Frattini subgroup, the Fitting subgroup, 
and a Sylow $p$-subgroup of $G$ respectively.
We use $G=[A]B$ to denote the semidirect product of $A$ 
and $B$, where $A$ is a normal subgroup of $G$. The set of prime divisors
of the order of $G$ is denoted by $\pi (G)$. As usual, 
$A_n$ and  $S_n$  are the alternating and the symmetric groups of
degree $n$ respectively. We use $E_{p^n}$ to denote an elementary 
abelian group of order $p^n$ and $Z_m$ to denote a cyclic group of
order $m$. Let $G$ be a
group of order $p_1^{a_1}p_2^{a_2} \ldots p_k^{a_k}$, where
$p_1>p_2> \ldots >p_k$. We say that $G$ has a Sylow tower
if there exists a series
$$
1=G_0\leq G_1\leq G_2\leq \ldots \leq G_{k-1}\leq G_k=G
$$ of normal subgroups of $G$ such that for each $i=1,2,\ldots
,k$, $G_{i}/G_{i-1}$ is isomorphic to a Sylow $p_i$-subgroup of
$G$. 
Recall that a positive integer $n$ is said to be squarefree if $n$
is not divisible by the square of any prime. 
A group is called metabelian if it contains
a normal abelian subgroup such that the corresponding quotent group
is also abelian.

We use $\frak H$ to denote the class of all groups $G$ such that 
each Schmidt subgroup $S$ of $G$ is a Hall subgroup of $G$.
It is clear that all nilpotent groups, all Schmidt groups, and all
squarefree groups belong to $\frak H$. If $T$ is biprimary 
non-nilpotent and  $T\in\frak H$, then $T$ is a Schmidt group.
Below the different examples of this class are given .

\medskip

{\bf Example 1.} Let $G=A\times B$, $(|A|, |B|)=1$, $A\in \frak H$,
$B \in \frak H$. Then, evidently, $G\in \frak H$.

\medskip

{\bf Example 2.}  Let $P$ be extraspecial of order
$409^3$. It is clear that $\Phi (P)=Z(P)=P^{\prime}$ 
has prime order $409$ and $P/\Phi (P)$
is elementary abelian of order $409^2$. The automorphism group of
$P/\Phi (P)$ is $GL(2,409)$. By theorem II.7.3 \cite{Hup}, $GL(2,409)$
has a cyclic subgroup $Z_{210}$ of order $5\cdot 41$.
Since $Z_{210}$ acts irreducibly on $P/\Phi(P)$, there
is $T=[P]Z_{210}$ such that $\Phi(p)=Z(T)$. 
The group $T$ possesses exactly three maximal subgroups:
$[P]Z_{5}$ is a Schmidt group; $[P]Z_{41}$ is a Schmidt group; $\Phi
(P)\times Z_{210}$ is a nilpotent subgroup. Therefore 
$\pi(T)=\{p,q,r\}$, where $p$, $q$, $r$ are distinct primes and
every Schmidt subgroup of $T$ is a Hall subgroup of $T$.

\medskip

{\bf Example 3.} Let
$G=[(E_4\times E_{25}\times E_7\times E_{121}\times E_{169}\times\cdots)]
Z_3$,
where $[E_4]Z_3$, $[E_{25}]Z_3$, $[E_7]Z_3$, $[E_{121}]Z_3$,
$[E_{169}]Z_3,\ldots$ are Schmidt groups in which all proper subgroups
are primary.
Let $K$ be a proper subgroup of $G$. If
$3$ does not divide $|K|$, then $K$ is nilpotent. Now suppose that $3$ 
divides $|K|$. Since $G$ is $p$-closed for any $p\ne 3$, it follows that
$K$ is 
$p$-closed too and there exists $[K_p]Z_3$. By Hall's theorem,
$[K_p]Z_3\subseteq [G_p]Z_3$. However all proper subgroups of
$[G_p]Z_3$ are primary. Thus $[K_p]Z_3=[G_p]Z_3$ is either Hall  
in $G$ or $K_p=1$. Since $p$ is an arbitrary prime number, $p\ne 3$, 
we see that $K$ is a Hall subgroup of $G$ and $G\in \frak H$.

\medskip

{\bf Lemma 1.}  {(\cite{Sm},\cite{MonUMK})}
{\it
Let $S$ be a Schmidt group. Then the following statements hold:

1) $S=[P]\langle y \rangle$, where $P$ is a normal Sylow $p$-subgroup,
$\langle y \rangle$~ is a non-normal cyclic Sylow $q$-subgroup,
$p$ and $q$ are distinct primes, $y^q\in Z(S)$;

2) $|P/P^{\prime}|=p^m$, where $m$ is the order of $p$ modulo $q$;

3) if $P$ is abelian, then $P$  is an elementary abelian $p$-group 
of order $p^m$ and $P$ is a minimal normal subgroup of $S$;  

4) if $P$ is non-abelian, then $Z(P)=P^{\prime}=\Phi (P)$ and
$|P/Z(P)| = p^m$;

5) if $P_1$ is a non-trivial normal $p$-subgroup of $S$ 
such that $P_1\ne P$,
then $P$ is non-abelian and $P_1 \subseteq Z(P)$;
 
6) $Z(S)=\Phi (S)=\Phi (P)\times \langle y^q\rangle$;
$S^{\prime}=P$, $P^{\prime}=(S^{\prime})^{\prime}=\Phi (P)$;

7) if $N$ is a proper normal subgroup of $S$, then
$N$ does not contain $\langle y \rangle$ and 
either $P\subseteq N$ or $N\subseteq \Phi (S)$.}

\medskip

We denote by $S_{\langle p,q\rangle}$-group a Schmidt group with
a normal Sylow $p$-subgroup and a cyclic Sylow $q$-subgroup. 

\medskip

{\bf Lemma 2.} (\cite{KnMon04}, lemma 2)
{\it If $K$ and $D$ are subgroups of $G$ such that 
$D$ is normal in $K$ and $K/D$ is an $S_{\langle p,q\rangle}$-subgroup, 
then each minimal supplement $L$ to $D$ in $K$ has the following
properties:

$1)$ $L$ is $p$-closed $\{p,q\}$-subgroup;

$2)$ all proper normal subgroups of $L$ are nilpotent;

$3)$ $L$ contains an $S_{\langle p,q\rangle}$-subgroup $[P]Q$ such that
$D$ does not contain $Q$ and $L=([P]Q)^L=Q^L$.}

\medskip

{\bf Lemma 3.} {\it If $G\in \frak H$, then every subgroup of $G$
and its every quotient group belongs to $\frak H$.}

\medskip

{\sc Proof.} Let $V\leq G\in \frak H$. 
If $V$ is non-nilpotent, then it contains a Schmidt subgroup $S$.
Since $G\in \frak H$, we can easily observe that $S$ is a Hall 
subgroup of $G$. It is clear that $S$ is a Hall subgroup of
$V$, hence $V\in \frak H$. 

Let $D$ be a normal subgroup of $G$ and $K/D$ is a Schmidt subgroup of
$G/D$. By the previous lemma, minimal supplement $L$ to $D$ in $K$ has
an $S_{\langle p,q\rangle}$-subgroup $[P]Q$ such that $D$ does not
include $Q$.
By lemma 1, $[P]QD/D$ is a Schmidt subgroup, hence
$[P]QD/D=K/D$. Since $G\in \frak H$, it follows that $[P]Q$ is a 
Hall subgroup of
$G$. Therefore $[P]QD/D=K/D$ is a Hall subgroup of $G/D$ and
$G/D\in \frak H$. The lemma is proved.

\medskip

{\bf Remark 1.} The class $\frak H$ is not closed
under direct products. 
For example, $S_3\in \frak H$, $Z_2\in \frak H$ but 
$S_3\times Z_2\not \in \frak H$. This shows that $\frak H$
is neither a formation nor a Fitting class.

\medskip

{\bf Lemma 4.}
{\it
1) If $G$ is not $p$-nilpotent, then $G$ has
a $p$-closed Schmidt $pd$-subgroup.

2) If $G$ is not $2$-closed, then $G$ has
a $2$-nilpotent Schmidt subgroup of even order.

3) If a $p$-solvable group $G$ is not $p$-closed, 
then $G$ has a $p$-nilpotent Schmidt $pd$-subgroup.}

\medskip

{\sc Proof.} 1. The proof of this part follows directly
from the Frobenius theorem (see, \cite{Hup}, theorem IV.5.4).

2. In \cite{Ber66} there is a proof based on Suzuki's theorem
of simple groups with independent Sylow $2$-subgroups.
Let us show another proof.
By induction, all proper subgroups of $G$ are 
$2$-closed. It follows that $G$ is not biprimary, 
(see part 1 of the lemma).  
If $G$ is solvable, then all biprimary Hall subgroups of 
$G$ are $2$-closed and $G$ is also $2$-closed, a contradiction.
Thus $G$ is not solvable. It is clear that
$G/\Phi(G)$ is a simple group.
Let $X$ be the conjugacy class of involutions of $G/\Phi(G)$.
By theorem IX.7.8 \cite{Hup2}, there exists involutions  
$x,y\in X$ such that $\langle x,y \rangle $ is not $2$-group.
It is well known that $\langle x,y \rangle $ is
the dihedral group of order $2|xy|$, (see \cite{Mon}, theorem 2.49).
It is not $2$-closed, a contradiction.

3. By theorem 5.3.13 \cite{Suz}, $G$ is a
$D_{\{p,q\}}$-group for any $q\in \pi (G)$. Suppose that $G$ is not
$p$-closed. Then $G$ contains a Hall $\{p,q\}$-subgroup $H$
such that $H$ is not $p$-closed for some prime $q\in \pi (G)$.
It is clear that $H$ is not $q$-nilpotent. By part 1 of the 
lemma, $H$ has a $p$-nilpotent $pd$-Schmidt subgroup. 
The lemma is proved.

\medskip

For any odd prime $p$ assertion 2 of lemma 3 is false.
If $p = 3$, then the counterexamples are $SL(2,2^n)$  
for any odd $n$ and $PSL(2,p)$ for $p \ge 5$.

\medskip

{\bf Lemma 5.} {\it If $G\in \frak H$, then $G$ 
possesses a Sylow tower.}


{\sc Proof.} First of all, we prove that if $G\in \frak H$ and $p$ is the
smallest prime dividing $|G|$, then $G$ is either $p$-closed or 
$p$-nilpotent.
Let $p=2$. If $G$ is not $2$-closed, then, by lemma 4 (2),
$G$ has a $2$-nilpotent Schmidt subgroup $S$ of even order.
Any Sylow $2$-subgroup of $S$ is cyclic. Since $G\in \frak H$, we deduce
that $S$ is a Hall subgroup of $G$ and Sylow $2$-subgroup of $G$ is cyclic.
Thus $G$ is $2$-nilpotent by theorem IV.2.8~\cite{Hup}.
Now suppose that $p>2$. Then $G$ is solvable. 
If $G$ is not $p$-closed, then, by lemma 3(3), $G$ 
has a $p$-nilpotent Schmidt $pd$-subgroup $T$. 
A Sylow $p$-subgroup $P$ of $T$ is cyclic. 
Since $G\in \frak H$, it follows that $T$ is a Hall subgroup of 
$G$ and $P$ 
is a Sylow $p$-subgroup of $G$. 
Thus $G$ is $p$-nilpotent by theorem IV.2.8~\cite{Hup}.

Therefore if $G\in \frak H$ and $p$ is the
smallest prime dividing $|G|$, then $G$ is either $p$-closed or 
$p$-nilpotent.
We use induction on $|G|$. Prove that $G$
possesses a Sylow tower.
Let $p$ be the smallest prime dividing $|G|$.
If a Sylow $p$-subgroup $P$ is normal in $G$, then, by lemma 3,
$G/P\in \frak H$ and, by induction, $G/P$
possesses a Sylow tower. 
Thus $G$ possesses a Sylow tower.
If $G$ is $p$-nilpotent, then
$G$ contains a normal subgroup $K$ such that $G/K$ is isomorphic
to a Sylow $p$-subgroup of $G$. By lemma 3, $K\in \frak H$ and, 
by induction, $K$ possesses a Sylow tower. 
Therefore $G$ possesses a Sylow tower.
The lemma is proved.

\medskip

{\bf Lemma 6.} {\it Let $G\in \frak H$ and $p$, $q$ are 
different prime divisors of $|G|$. 
Then any Hall $\{p,q\}$-subgroup of $G$ is either nilpotent
or Schmidt group.}

\medskip

{\sc Proof.} By lemma 5, $G$ is solvable, so $G$ has a
Hall $\{p,q\}$-subgroup $K$. Assume that $K$ is non-nilpotent. 
Then $K$ contains a Schmidt subgroup $S$. Since $G\in \frak H$, it
implies that $S$ must be a Hall subgroup of $G$. Therefore $S=K$.
The lemma is proved.

\medskip

{\bf Lemma 7.} {\it Let $n\ge 2$ be a positive integer
and $p$ be a prime number. 
Denote by $\pi$ the set of prime numbers $q$ such that $q$ divides 
$p^n-1$, but $q$ does not divide $p^{n_1}-1$ for all $1\leq n_1<n$. 
Then $GL(n,p)$ has a cyclic Hall $\pi$-subgroup.}

\medskip

{\sc Proof.} The group $G=GL(n,p)$ has order
$$
p^{n(n-1)/2}(p^n-1)(p^{n-1}-1)...(p^2-1)(p-1).
$$
By theorem II.7.3 \cite{Hup}, $G$ contains a cyclic subgroup
$T$ of order $p^n-1$. Denote by $T_\pi$ a Hall $\pi$-subgroup of $T$.
Since $q$ does not divide $p^{n_1}-1$ for all $q\in\pi$ 
and all $1\leq n_1<n$, it follows that $T_\pi$ is a 
Hall $\pi$-$\pi$-subgroup of $G$.
The lemma is proved.


\section{Proof of Theorem and Corollary}

{\bf Theorem.} {\sl Let $G$ be a finite non-nilpotent group
in which every Schmidt subgroup is Hall in $G$.
Then the following statements hold:

1) if $P$ is a non-normal Sylow $p$-subgroup of $G$, then
$P$ is cyclic and every maximal subgroup of $P$ is contained in $Z(G)$;

2) if $P$ is a normal Sylow $p$-subgroup of $G$ and $G$ is 
not $p$-decomposable, then $P$ is either minimal normal in $G$  
or nonabelian, $Z(P)=P^{\prime}=\Phi (P)$,
and $P/\Phi (P)$ is minimal normal in $G/\Phi (P)$;

3) if $P_1$ is a normal $p$-subgroup of $G$, $P_1$ other than a 
Sylow one and $G$ is not $p$-decomposable, then $P_1$ is
contained in $Z(G)$;

4) if $Z(G)=1$, then $G$ has a normal abelian Hall subgroup 
$A$ in which every Sylow subgroup is minimal normal 
in $G$, $G/A$ is cyclic and $|G/A|$ is squarefree.}

\medskip

{\sc Proof.}
1. Let $G\in\frak H$ and $p\in\pi(G)$. Assume that $G$ has 
a non-normal Sylow $p$-subgroup $P$. By lemma 5, $G$ is solvable,
hence $G$ contains a Hall $\{p,q\}$-subgroup for any 
$q\in \pi(G)\setminus \{p\}$ by theorem 5.3.13 \cite{Suz}.
Since $P$ is non-normal in $G$, it follows that $G$ contains a 
not $p$-closed Hall $\{p,q\}$-subgroup $K$ for some 
$q\in\pi(G)\setminus \{p\}$. By lemma 4, $K$ has a $q$-closed
Schmidt subgroup $Q$. Under the condition of $G\in\frak H$,
$Q$ is the same as $K$.  
By the properties of Schmidt groups (see lemma 1(1)),
every Sylow $p$-subgroup of $K$ is cyclic.
Since $K$ is a Hall subgroup of $G$, we see that a Sylow 
$p$-subgroup of $K$ is a Sylow subgroup of $G$. Thus $P$ is cyclic.

Let $P_1$ be a maximal subgroup of $P$. If $P_1=1$, then $P_1\subseteq
Z(G)$. Assume that  $P_1\ne 1$. It is clear that $G$ has a
Hall $\{p,q\}$-subgroup $PQ$ for any prime $q\in \pi (G)\setminus \{p\}$,
where $Q$ is some Sylow $q$-subgroup of $G$. If $PQ$ is nilpotent,
then $Q\subseteq C_G(P_1)$. If $PQ$ is non-nilpotent, then  
$PQ$ is a Schmidt group by lemma 6.
If $PQ$ is $p$-closed, then $P$ has a prime order by lemma 1(3), 
a contradiction.
Hence $PQ$ is $q$-closed and $P_1\subseteq  Z(PQ)$ by lemma 1(1),
i.e. $Q\subseteq  C_G(P_1)$. Thus $C_G(P_1)$ contains a Sylow
$q$-subgroup for every $q\in \pi (G)\setminus \{p\}$. Since
$P\subseteq  C_G(P_1)$, we have $C_G(P_1)=G$ and $P_1\subseteq  Z(G)$.

2. Let Sylow $p$-subgroup $P$ be a normal subgroup of $G$.
Suppose that $P$ is not a minimal normal subgroup of $G$. 
In particular, $|P|>p$. By Schur-Zassenhaus theorem, $G$ has a
Hall $p^{\prime}$-subgroup $H$. By the hypothesis of the theorem, 
$G$ is not $p$-decomposable. Hence $H$ has a Sylow subgroup $Q$
such that $[P]Q$ is non-nilpotent. By lemma 6, $[P]Q$ is a Schmidt 
subgroup. By our assumption, $P$ is not minimal normal in $G$, 
it follows that $P$ is not minimal normal in $[P]Q$.
By the properties of Schmidt groups (see lemma 1(3)),
$P$ is non-abelian  and $Z(P)=P^{\prime}=\Phi (P)$.
Since $[P/\Phi (P)](Q\Phi (P)/\Phi (P)$ is a Schmidt group,
$P/\Phi (P)$ is its minimal normal subgroup. We see that  
$P/\Phi(P)$ is a minimal normal subgroup of $G/\Phi(P)$. 
The statement 2 is proved.

3. We denote by $G_p$ a Sylow $p$-subgroup of $G$. Assume that $Z(G)$ 
does not contain $P_1$. Then $|P_1|\geq p$, $|G_p|\geq p^2$, and $G_p$ 
is normal in $G$ by claim 1 of the theorem.
Let $G_q$ be a Sylow $q$-subgroup of $G$,
$q\in \pi (G)\setminus \{p\}$. By lemma 6, the product $G_pG_q$ either 
nilpotent or a Schmidt group. Suppose $G_pG_q$ is nilpotent
for all $q\in \pi(G)\setminus \{p\}$. In this case, 
$G=G_p\times G_{p^\prime}$, a contradiction. Thus our assumption
is false and there exist a prime $r\in \pi (G)\setminus \{p\}$ 
such that $G_pG_r$ is non-nilpotent. It follows that
$G_pG_r$ is a $p$-closed Schmidt group and $P_1$ is its
normal $p$-subgroup. 
By the properties of Schmidt groups (see lemma 1(3)),
$P_1\subseteq Z(G_pG_r)$. 
Thus, $P_1\subseteq Z(G_p)$ and $G_r\subseteq
C_G(P_1)$ for all $r\in \pi (G)\setminus \{p\}$ such that
$G_pG_r$ is not nilpotent.
If $G_pG_r$ is nilpotent, then $G_q\subseteq
C_G(P_1)$. Therefore $P_1\subseteq Z(G)$.

4. We denote by $\frak A$, $\frak N$ and $\frak E$ the classes of
all abelian, all nilpotent, and all finite groups respectively.
We define $\frak N\circ\frak A=\{G\in \frak E\mid G^{\frak A}\in \frak N\}$
and call $\frak N\circ\frak A$ the product of classes $\frak N$ and 
$\frak A$, where $G^{\frak A}$ denotes $\frak A$-residual of $G$, 
i.e. the smallest normal subgroup of $G$ quotient by which belogs 
to $\frak A$. 
The other definitions and terminology
about formations could be referred to Doerk, Hawkes (1992),
Huppert (1967) and Shemetkov (1978).
It is clear that $G^{\frak A}=G^{\prime }$ is the derived subgroup of $G$.
Hence $\frak N\circ\frak A$ consists of all groups $G$ whose the
derived groups are nilpotent. The class $\frak N\circ\frak A$ 
is a saturated formation.
Now, by induction on $|G|$, we prove that  
$\frak H\subseteq\frak N\circ\frak A$. Suppose the assertion is false.
Let $G$ be a counterexample of minimal order and  
$G\in\frak H\setminus \frak N\circ\frak A$. 
By lemma 5, $G$ is solvable and, by lemma 3, $G/N\in \frak H$ 
for every normal subgroup $N\ne 1$ of $G$. By induction, 
$G/N\in \frak N\circ\frak A$.
Since $\frak N\circ\frak A$ is a saturated formation, it follows that $G$
is primitive (see (\cite{Mon}, p. 143).  By theorem 4.42 \cite{Mon}, 
$F=F(G)=C_G(F)\simeq E_{p^n}$ is a minimal normal subgroup of $G$
and, by the above claim 3 of the theorem, $F$ is a Sylow subgroup of $G$.

If $n=1$, then $G/F$ is isomorphic to a subgroup of the automorphism
group of $F$, where $|F|=p$. Thus $G\in \frak N\circ\frak A$. 
Next, we assume that $n\geq 2$. Since $[F]G_q$ is a Hall non-nilpotent
subgroup of $G$, we have, by lemma 6, that $[F]G_q$ is a 
Schmidt subgroup for every $q\in \pi =\pi(G/F)$. By lemma 1(2), 
$q$ divides $p^n-1$, but $q$ does not divide $p^{n_1}-1$ for 
all $1\leq n_1<n$. 
The quotient group $G/F$ is isomorphic to a subgroup $K$ of $GL(n,p)$, 
$K$ has a cyclic Hall $\pi$-subgroup $T$ by lemma 7. 
By theorem 5.3.2 \cite{Suz}, $G/F$ is contained in some subgroup
$T^x$, $x\in GL(n,p)$. Thus $G/F$ is a cyclic and $G\in \frak N\circ\frak A$.
Let $G\in \frak H$ and $Z(G)=1$. Then $G$ is not $p$-decomposable
for any $p\in \pi (G)$. The assertion (3) implies that every 
minimal normal subgroup of $G$ is a Sylow subgroup of $G$. So
$F(G)=A$ is an abelian Hall subgroup of $G$ in which 
every Sylow subgroup is minimal normal in $G$.

Let $B$ be a complement to $A$ in $G$. By 
Schur-Zassenhaus theorem, (see \cite{Mon}, p. 136), in any case, $G$  
has some subgroup $B$ such that $G=AB$ and $A\cap B=1$.
The assertion 1 implies that $|B|$ is squarefree.
Since $G\in \frak N\circ\frak A$, it follows that $B$ is abelian. 
Therefore $B$ is cyclic. The theorem is proved.
 
\medskip

{\bf Corollary. } {\sl Let $G$ be a finite non-nilpotent group
in which every Schmidt subgroup is a Hall subgroup of $G$.
Then $G$ contains a nilpotent Hall subgroup $H$ such that 
$G/H$ is cyclic. In particular, $G/\Phi (G)$ is metabelian.}

\medskip

{\sc Proof.} If $Z(G)=1$, then the claim of the corollary
is the same as assertion 4 of the theorem. 
Let $Z(G)\ne 1$. Denote by $N$ a subgroup of prime order 
$p$, $N\subseteq Z(G)$. By induction, we have 
$\overline G=[A/N](B/N)$, where $A/N$ is a nilpotent Hall
subgroup of $G/N$ and $B/N$ is cyclic.
Since $N\subseteq Z(G)$, we see that $A$ and $B$ are nilpotent,
(see \cite{Mon}, lemma 3.15). If $A$ is a Hall subgroup of $G$, then, by
Schur-Zassenhaus theorem, $B=N\times B_1$ and $G=[A]B_1$, 
where $A$ is a nilpotent Hall subgroup of $G$ and $B_1$ is a cyclic
subgroup. In this case, the corollary is prove.
Now we assume that $A$ is not Hall subgroup of $G$. Then $A=N\times  A_1$, 
where $A_1$ is a normal nilpotent Hall subgroup of $G$ and
$G=[A_1]B$. Denote by $B_1$ the product of all Sylow subgroups
$B_{p_i}$ of $B$ such that $B_{p_i}$ are normal in $G$ for all $i$.
Respectively, denote by $B_2$ the product of all Sylow subgroups
$B_{p_i}$ of $B$ such that $B_{p_i}$ are non-normal in $G$ for all $i$.
It is clear that  
$G=[A\times B_1]B_2$, where $A\times B_1$ is a normal Hall subgroup
of $G$ and all Sylow subgroups of $B_2$ are cyclic by 
the assertion 1 of the theorem.
Since $B_2$ is nilpotent, it follows that $B_2$ is cyclic.
Therefore in any case, $G$ contains a nilpotent Hall subgroup $H$
such that $G/H$ is cyclic. Since $\Phi (H)\subseteq \Phi (G)$,
$H/\Phi (H)$ is abelian, we see that $G/\Phi (G)$ is metabelian.
The corollary is proved.

\medskip

{\bf Remark 2.} For any natural number $n\ge 3$ there exists a
nilpotent subgroup $A$ such that the nilpotent length of $A$ is 
equal to $n$. Let $p$ and $q$
are  distinct primes and $p,q\notin\pi(A)$. 
By theorem 1.3 \cite{MonUMK}, there exists an
$S_{\langle p,q\rangle}$-subgroup $B$. 
All Schmidt subgroups of $G=A\times B$ are Hall subgroups of $G$ and 
the derived length of $G$ is equal to $n$. Now, if $G$ is a non-nilpotent 
group and $G\in\frak H$, then its derived length is not bounded above.

\bigskip


\noindent V.\,N. KNIAHINA

\noindent  Department of mathematics, 
Gomel Engineering Institute, Gomel
246035, BELARUS

\noindent E-mail address: knyagina@inbox.ru

\bigskip

\noindent V.\,S. MONAKHOV

\noindent Department of mathematics, Gomel F. Scorina State
University, Gomel 246019,  BELARUS 

\noindent E-mail address: Victor.Monakhov@gmail.com

\end{document}